\documentclass[12pt,twoside,reqno]{amsart}
\usepackage[colorlinks=true,citecolor=blue]{hyperref}
\usepackage{mathptmx, amsmath, amssymb, amsfonts, amsthm, mathptmx, enumerate, color}
\usepackage{graphicx}
\usepackage{epstopdf}

\newtheorem{theorem}{Theorem}[section]

\newtheorem{lemma}{Lemma}[section]

\theoremstyle{definition}

\newtheorem{remark}{Remark}[section]

\numberwithin{equation}{section}

\begin{document}
\setcounter{page}{1}

\vspace*{1.0cm}
\title[Norm of the Generalized Hilbert Matrix Operator Between Some Spaces of Analytic Functions]
{Norm of the Generalized Hilbert Matrix Operator Between Some Spaces of Analytic Functions}
\author[S. Feng, S. Ye]{ Shaoyang Feng, Shanli Ye$^{*}$}
\maketitle
\vspace*{-0.6cm}

\begin{center}
{\footnotesize {\it

School of Science, Zhejiang University of Science and Technology, Hangzhou 310023, China.

}}\end{center}

\vskip 4mm {\small\noindent {\bf Abstract.}
We study the generalized Hilbert operator
\[
\mathcal H_bf(z)=\int_0^1\frac{(1-t)^b}{(1-tz)^{b+1}}f(t)\,dt,
\qquad b>0,
\]
on several spaces of analytic functions in the unit disk.
For logarithmically weighted Korenblum spaces, we obtain an exact
supremum representation and a lower estimate for the norm of the
mapping from $H_{\alpha,\log}^{\infty}$ to $H_\alpha^\infty$.
We then establish two-sided estimates and characterize the boundedness
range of $\mathcal H_b$ on the $\alpha$-Bloch spaces $\mathcal B^\alpha$.
We also derive norm bounds for the mapping from $H^\infty$ into the
endpoint Zygmund-type space $\Lambda_1^{1,*}$.
Furthermore, we give a norm representation for
$\mathcal H_b:H_\alpha^\infty\to\mathcal B^{\alpha+1}$ and identify
a norm-attaining function. In particular, for $0<\alpha\leq2/3$,
the norm is expressed explicitly in terms of a Beta function.
Finally, we obtain an exact supremum representation of the norm of
$\mathcal H_b$ on $H_{\alpha,\log}^{\infty}$, together with explicit
norm estimates.  Our results extend those of Hao et al. [J. Geom. Anal. 35 (2025), Paper No. 184]
for the classical Hilbert matrix operator ($b=0$) to the generalized case $b>0$.

\noindent {\bf Keywords.}
Generalized Hilbert operator, operator norm, Korenblum space, Bloch-type space, Zygmund-type space }

\renewcommand{\thefootnote}{}
\footnotetext{ $^*$Corresponding author.
\par
E-mail addresses:  slye@zust.edu.cn (S. Ye).
\par
}

\section{Introduction}
 The Hilbert matrix, which we denote by $\mathcal H$ when regarded as an operator, is a classical example of an infinite matrix whose
properties depend on the space on which it acts. Magnus \cite{Magnus1950}
studied its spectrum on $\ell^2$, and the spectral properties of
related Hilbert matrices have subsequently been examined in other
settings; see, for example, \cite{AlemanMontesSarafoleanu2012}.
The action of the Hilbert matrix on analytic function spaces was studied
by Diamantopoulos and Siskakis \cite{DiamantopoulosSiskakis2000}, who
proved its boundedness on the Hardy spaces $H^p$ for $1<p<\infty$ and
obtained the estimate
\[
\|\mathcal H\|_{H^p\to H^p}\leq \frac{\pi}{\sin(\pi/p)},
\qquad 2\leq p<\infty.
\]
The exact norm on $H^p$, $1<p<\infty$, was later established by
Dostani\'c, Jevti\'c and Vukoti\'c \cite{DostanicJevticVukotic2008}.
On the Bergman spaces $A^p$, Diamantopoulos
\cite{Diamantopoulos2004} obtained upper estimates, including
\[
\|\mathcal H\|_{A^p\to A^p}\leq \frac{\pi}{\sin(2\pi/p)},
\qquad 4\leq p<\infty.
\]
Dostani\'c, Jevti\'c and Vukoti\'c \cite{DostanicJevticVukotic2008}
proved that this bound is sharp when $4\leq p<\infty$ and conjectured
the same formula for $2<p<4$. Bo\v zin and Karapetrovi\'c
\cite{BozinKarapetrovic2018} settled that conjecture, giving
\[
\|\mathcal H\|_{A^p\to A^p}
=\frac{\pi}{\sin(2\pi/p)},\qquad 2<p<\infty.
\]

The norm problem has also been considered on spaces with pointwise
growth conditions. Lindstr\"om, Miihkinen and Wikman
\cite{LindstromMiihkinenWikman2019} determined the exact norm of
$\mathcal H$ on the Korenblum space $H_\alpha^\infty$ for
$0<\alpha\leq 2/3$ and obtained an upper bound for $2/3<\alpha<1$.
Dai \cite{Dai2022} subsequently studied the norm on the remaining range
and identified a transition beyond which the previous candidate for
the exact norm ceases to be sharp. Diamantopoulos and Siskakis
\cite{DiamantopoulosSiskakis2000} also showed that $\mathcal H$ is not
bounded on $H^\infty$. Nevertheless, the operator admits bounded
mappings into suitable spaces of analytic functions. Hu and Ye
\cite{HuYe2025} obtained norm results for logarithmically weighted
Korenblum spaces, including the mapping into $H_\alpha^\infty$ and the
operator acting on $H_{\alpha,\log}^\infty$. They also studied the
operator on the $\alpha$-Bloch spaces. In a related direction, Guo and
Tang \cite{GuoTang2024} studied the mapping from $H^\infty$ into an
endpoint Zygmund-type space and established two-sided norm estimates.
Hu and Ye \cite{HuYe2025} also considered the mapping from
$H_\alpha^\infty$ into $\mathcal B^{\alpha+1}$. More recently, Ye and Zheng \cite{ysl} introduced the logarithmically weighted Bloch and Hardy
spaces relevant to the present work and studied the Hilbert matrix operator on these spaces.

Let $\mathbb D=\{z\in\mathbb C:|z|<1\}$ and let $H(\mathbb D)$ denote
the space of analytic functions on $\mathbb D$. For
$f(z)=\sum_{k=0}^\infty a_kz^k\in H(\mathbb D)$, the Hilbert matrix
operator is formally defined by
\[
\mathcal H f(z)=\sum_{n=0}^\infty
\left(\sum_{k=0}^\infty\frac{a_k}{n+k+1}\right)z^n,
\]
whenever the indicated series are convergent. Equivalently, on
function classes for which the integral is well defined, it has the
representation \cite{DiamantopoulosSiskakis2000}
\[
\mathcal H f(z)=\int_0^1\frac{f(t)}{1-tz}\,dt,\qquad z\in\mathbb D.
\]
For example, the Fej\'er--Riesz inequality guarantees the required
integrability when $f\in H^1$; see \cite{Duren1970}.

We consider the generalized Hilbert operator
\begin{equation}\label{eq:Hb-introduction}
\mathcal H_bf(z)=\int_0^1\frac{(1-t)^b}{(1-tz)^{b+1}}f(t)\,dt,
\qquad b>0,\quad z\in\mathbb D,
\end{equation}
whenever the integral converges. The generalized Hilbert operator was
introduced in a broader setting by Li and Stevi\'c \cite{LiStevic2009}.
For $b=0$, the expression in \eqref{eq:Hb-introduction} reduces to the
classical Hilbert matrix operator. Recently, Li, Pan and Zhou
\cite{LiPanZhou2026} determined the norm of $\mathcal H_b$ on Hardy
spaces and, in particular, proved
\[
\|\mathcal H_b\|_{H^\infty\to\mathcal B}
=\frac1{b+1}+2,\qquad b>0.
\]

These results lead us to ask how the norm estimates for the classical
Hilbert matrix operator change when it is replaced by $\mathcal H_b$.
The factor $(1-t)^b$ alters the behavior of the integral kernel near
$t=1$, so both the boundedness ranges and the constants in the estimates
may depend on $b$. We first use a weighted composition representation
to study $\mathcal H_b$ as an operator from the logarithmically weighted
Korenblum space $H_{\alpha,\log}^\infty$ into $H_\alpha^\infty$.
We next investigate its action on the $\alpha$-Bloch spaces, where
growth estimates and test functions yield boundedness criteria and
two-sided norm estimates. We then consider the mapping from
$H^\infty$ to the endpoint Zygmund-type space $\Lambda_1^{1,*}$,
whose norm involves integral means of the second derivative, followed
by the mapping from $H_\alpha^\infty$ into $\mathcal B^{\alpha+1}$,
where the norm estimates are closely related to Beta integrals.
Finally, we study the norm of $\mathcal H_b$ acting on
$H_{\alpha,\log}^\infty$ itself.

We study the action of $\mathcal H_b$ between several spaces of
analytic functions. The paper is organized as follows. In Section~2 we
recall the function spaces and notation used below. Section~3 studies
the mapping $\mathcal H_b:H_{\alpha,\log}^\infty\to H_\alpha^\infty$,
including its weighted composition representation and norm estimates.
Section~4 is devoted to boundedness and norm estimates on the
$\alpha$-Bloch spaces $\mathcal B^\alpha$. In Section~5 we obtain
two-sided norm bounds for $\mathcal H_b:H^\infty\to\Lambda_1^{1,*}$.
Section~6 considers the mapping
$\mathcal H_b:H_\alpha^\infty\to\mathcal B^{\alpha+1}$, including a
norm representation and estimates for the corresponding operator
norm. Finally, Section~7 is concerned with the norm of $\mathcal H_b$
on $H_{\alpha,\log}^\infty$.

 \section{Preliminaries}

Let $\mathbb D$ denote the open unit disk of the complex plane
$\mathbb C$, and let $H(\mathbb D)$ denote the set of all analytic
functions in $\mathbb D$.

For $0<p\leq\infty$, the Hardy space $H^p$ consists of all
functions $f\in H(\mathbb D)$ such that
\[
\|f\|_{H^p}=\sup_{0\leq r<1}M_p(r,f)<\infty,
\]
where
\[
M_p(r,f)=\left(\frac{1}{2\pi}\int_0^{2\pi}
|f(re^{it})|^p\,dt\right)^{1/p},\qquad 0<p<\infty,
\]
and
\[
M_\infty(r,f)=\sup_{0\leq t<2\pi}|f(re^{it})|.
\]
We refer to \cite{Duren1970} for the terminology and results
on Hardy spaces. In particular, $H^\infty$ denotes the space
of bounded analytic functions on $\mathbb D$.

For $0<\alpha<1$, the Korenblum space $H_\alpha^\infty$
consists of all functions $f\in H(\mathbb D)$ satisfying
\[
\|f\|_{H_\alpha^\infty}
=\sup_{z\in\mathbb D}(1-|z|^2)^\alpha|f(z)|<\infty.
\]

 From \cite{HuYe2025}, we recall that  the logarithmically weighted Korenblum space
$H_{\alpha,\log}^{\infty}$ consists of those
$f\in H(\mathbb D)$ such that
\[
\|f\|_{H_{\alpha,\log}^\infty}=
\sup_{z\in\mathbb D}(1-|z|^2)^\alpha
\log\frac{2e^{1/\alpha}}{1-|z|^2}|f(z)|<\infty.
\]
It is easy to see that $H^\infty\subsetneq H_{\alpha,\log}^\infty
\subsetneq H_\alpha^\infty$.

For $0<\alpha<\infty$, the $\alpha$-Bloch space $\mathcal B^\alpha$
consists of all functions $f\in H(\mathbb D)$ with
\[
\|f\|_{\mathcal B^\alpha}
=|f(0)|+\sup_{z\in\mathbb D}(1-|z|^2)^\alpha|f'(z)|<\infty.
\]
The case $\alpha=1$ is the classical Bloch space
$\mathcal B=\mathcal B^1$; see \cite{AndersonCluniePommerenke1974,Zhu1993}.

Let $1\leq p<\infty$ and $0<\alpha\leq1$. The mean Lipschitz
space $\Lambda_\alpha^p$ consists of those functions $f\in H(\mathbb D)$
having  non-tangential limits almost everywhere such that
$\omega_p(t,f)=O(t^\alpha)$ as $t\to0$. Here $\omega_p(\cdot,f)$ denotes the integral modulus of continuity of order
 $p$ of $f(e^{i\theta})$. It is known (see \cite{Duren1970}) that
$\Lambda_\alpha^p$ is a subset of $H^p$ and
\[
\Lambda_\alpha^p
=\left\{f\in H(\mathbb D):
M_p(r,f')=O\!\left(\frac{1}{(1-r^2)^{1-\alpha}}\right),
\text{ as }r\to1^-\right\}.
\]
The space $\Lambda_\alpha^p$ is a Banach space with the norm
$\|\cdot\|_{\Lambda_\alpha^p}$ given by
\[
\|f\|_{\Lambda_\alpha^p}
=|f(0)|+\sup_{0\leq r<1}(1-r^2)^{1-\alpha}M_p(r,f').
\]

For $0<p\leq\infty$, the Zygmund-type space $\mathcal Z_p$
is the space of $f\in H(\mathbb D)$ such that
\[
\|f\|_{\mathcal Z_p}
=|f(0)|+|f'(0)|+
\sup_{0\leq r<1}(1-r^2)M_p(r,f'')<\infty.
\]
When $p=\infty$, the space $\mathcal Z_\infty$ is the well-known
Zygmund space (see \cite{Stevic2009}). The space $\mathcal Z_1$, on the other hand, is closely related
to the mean Lipschitz space $\Lambda_{1/p}^p$.

\section[Norm of the generalized Hilbert matrix operator]{Norm of the generalized Hilbert matrix operator
$\|\mathcal H_b\|_{H_{\alpha,\log}^{\infty}\to H_{\alpha}^{\infty}}$}
\label{sec:Hb-representation}

For $b>0$ and $0<\alpha<1$, since
$H_{\alpha,\log}^{\infty}\subseteq H_{\alpha}^{\infty}$,  it is natural to consider
$\mathcal H_b$ as an operator from the logarithmically weighted Korenblum
space $H_{\alpha,\log}^{\infty}$ into the Korenblum space
$H_{\alpha}^{\infty}$.  It is known that the classical Hilbert matrix operator
$\mathcal H$ is bounded on $H_{\alpha}^{\infty}$ for $0<\alpha<1$.
Moreover, for $0<t<1$ and $z\in\mathbb D$,
\[
\left|\frac{(1-t)^b}{(1-tz)^{b+1}}\right|
=\frac{1}{|1-tz|}\left(\frac{1-t}{|1-tz|}\right)^b
\leq\frac{1}{1-t|z|}.
\]
Consequently, by comparison with the classical kernel,
$\mathcal H_b$ is also bounded from
$H_{\alpha,\log}^{\infty}$ into $H_{\alpha}^{\infty}$.
In this section, we aim to derive norm estimates for
$\mathcal H_b:H_{\alpha,\log}^{\infty}\to H_{\alpha}^{\infty}$
when $b>0$ and $0<\alpha<1$.

Let $f(z)=\sum_{k=0}^{\infty}a_kz^k\in H_{\alpha,\log}^{\infty}$.
For $b>0$,  the matrix associated with the generalized Hilbert operator $\mathcal H_b$
 has entries
\[
\frac{\Gamma(n+b+1)\Gamma(n+k+1)}
{\Gamma(n+1)\Gamma(n+k+b+2)},\qquad n,k\geq0.
\]
Its action on $f$ admits  the following integral representation:
\begin{align}
\mathcal H_bf(z)
&=\sum_{n=0}^{\infty}\left(
\sum_{k=0}^{\infty}
\frac{\Gamma(n+b+1)\Gamma(n+k+1)}
{\Gamma(n+1)\Gamma(n+k+b+2)}a_k
\right)z^n \notag\\
&=\sum_{n=0}^{\infty}\left(
\frac{(b+1)_n}{n!}
\int_0^1t^n(1-t)^bf(t)\,dt
\right)z^n \notag\\
&=\int_0^1(1-t)^bf(t)
\sum_{n=0}^{\infty}\frac{(b+1)_n}{n!}(tz)^n\,dt \notag\\
&=\int_0^1\frac{(1-t)^b}{(1-tz)^{b+1}}f(t)\,dt,
\qquad z\in\mathbb D.
\label{eq:Hb-integral}
\end{align}
Here $(b+1)_n$ denotes the Pochhammer symbol. For fixed
$z\in\mathbb D$, the interchange of the $n$-series and the integral
follows from the binomial expansion and the integrable majorant
$(1-t)^{b-\alpha}$. The inner $k$-series can be identified with the
moment integral by first replacing $f(t)$ with $f(\rho t)$,
$0<\rho<1$, and then letting $\rho$ tend to $1$ from below. In fact, it converges
in the usual sense for the present parameters: the standard growth
estimate for Taylor partial sums of $f\in H_\alpha^\infty$ gives
$\sum_{k=0}^{N}a_k=O((N+1)^\alpha\log(N+2))$, while
$B(n+k+1,b+1)=O((k+1)^{-b-1})$ and its consecutive difference is
$O((k+1)^{-b-2})$. Summation by parts applies because
$\alpha<b+1$.

We also have the weighted composition representation
\begin{equation}
\mathcal H_bf(z)=\int_0^1T_tf(z)\,dt,
\qquad z\in\mathbb D,
\label{eq:Hb-weighted-composition}
\end{equation}
where
\[
T_tf(z)=w_t(z)f\bigl(\phi_t(z)\bigr),\qquad
w_t(z)=\frac{(1-t)^b}{1-(1-t)z},\qquad
\phi_t(z)=\frac{t}{1-(1-t)z}.
\]
Indeed, for real $z=r\in[0,1)$, the substitution
$x=\phi_t(r)$ gives
\[
t=\frac{x(1-r)}{1-xr},\qquad
1-t=\frac{1-x}{1-xr},\qquad
dt=\frac{1-r}{(1-xr)^2}\,dx,
\]
and hence
\[
\int_0^1w_t(r)f\bigl(\phi_t(r)\bigr)\,dt
=\int_0^1\frac{(1-x)^b}{(1-xr)^{b+1}}f(x)\,dx.
\]
Both sides of \eqref{eq:Hb-weighted-composition} are analytic in
$z\in\mathbb D$; the identity theorem therefore extends the equality
from $[0,1)$ to the whole unit disk. Notice that setting $b=0$ in
\eqref{eq:Hb-integral}--\eqref{eq:Hb-weighted-composition} recovers
the classical Hilbert matrix operator and its weighted composition
representation.

\begin{theorem}\label{thm:Hb-radial-norm}
For $b>0$ and $0<\alpha<1$, we have
\[
\|\mathcal H_b\|_{H_{\alpha,\log}^{\infty}\to H_\alpha^{\infty}}
=\sup_{0\leq r<1}\int_0^1
\frac{(1-t)^{b-\alpha}(1+r)^\alpha[1-(1-t)r]^{2\alpha-1}}
{[1+t-(1-t)r]^\alpha\log\!\frac{2e^{1/\alpha}}{1-\phi_t(r)^2}}\,dt,
\]
and
\[
\|\mathcal H_b\|_{H_{\alpha,\log}^{\infty}\to H_\alpha^{\infty}}
\geq\int_0^1
\frac{(1-t)^b}
{(1-t^2)^\alpha\log\!\frac{2e^{1/\alpha}}{1-t^2}}\,dt.
\]
\end{theorem}

\begin{proof}
Let $b>0$ and $0<\alpha<1$. We use the test function
\[
f_\alpha(z)=\frac{1}
{(1-z^2)^\alpha\log\!\frac{2e^{1/\alpha}}{1-z^2}}.
\]
The analytic branches are chosen to agree with the positive real
values on $(-1,1)$. Hu and Ye \cite{HuYe2025} established that this function
belongs to $H_{\alpha,\log}^{\infty}$ and satisfies
\[
\|f_\alpha\|_{H_{\alpha,\log}^{\infty}}=1.
\]

The weighted composition operator $T_t$ applied to the function
$f_\alpha$ can be written as
\[
\begin{aligned}
T_tf_\alpha(z)
&=\frac{(1-t)^b}{1-(1-t)z}
 f_\alpha\!\left(\frac{t}{1-(1-t)z}\right)\\
&=\frac{(1-t)^b}{1-(1-t)z}
\frac1{\left(1-\left(\frac{t}{1-(1-t)z}\right)^2\right)^\alpha
\log\!\frac{2e^{1/\alpha}}
{1-\left(\frac{t}{1-(1-t)z}\right)^2}}.
\end{aligned}
\]
Since $\mathcal H_bf_\alpha(z)=\displaystyle\int_0^1T_tf_\alpha(z)\,dt$
and $\|f_\alpha\|_{H_{\alpha,\log}^{\infty}}=1$, we get that
\begin{align}
&\|\mathcal H_b\|_{H_{\alpha,\log}^{\infty}\to H_\alpha^{\infty}}
\geq\|\mathcal H_bf_\alpha\|_{H_{\alpha}^{\infty}}
=\sup_{z\in\mathbb D}(1-|z|^2)^\alpha
\left|\int_0^1T_tf_\alpha(z)\,dt\right|\notag\\
&\geq\sup_{0\leq r<1}\int_0^1
(1-r^2)^\alpha\frac{(1-t)^b}{1-(1-t)r}
\frac{dt}
{\left(1-\left(\frac{t}{1-(1-t)r}\right)^2\right)^\alpha
\log\!\frac{2e^{1/\alpha}}
{1-\left(\frac{t}{1-(1-t)r}\right)^2}}\notag\\
&=\sup_{0\leq r<1}\int_0^1
\frac{(1-t)^b(1+r)^\alpha\bigl((t-1)r+1\bigr)^{2\alpha-1}}
{(1-t)^\alpha\bigl((t-1)r+1+t\bigr)^\alpha
\log\!\frac{2e^{1/\alpha}}
{1-\left(\frac{t}{1-(1-t)r}\right)^2}}\,dt.
\label{eq:Hb-lower}
\end{align}
The last equality follows from
\[
1-\left(\frac{t}{1-(1-t)r}\right)^2
=\frac{(1-r)(1-t)\bigl((t-1)r+1+t\bigr)}
{\bigl((t-1)r+1\bigr)^2}.
\]

On the other hand, we have that
\begin{align}
|T_tf(z)|
&=|w_t(z)f(\phi_t(z))|\notag\\
&=|w_t(z)|\,
\frac1{(1-|\phi_t(z)|^2)^\alpha
\log\!\frac{2e^{1/\alpha}}{1-|\phi_t(z)|^2}}
(1-|\phi_t(z)|^2)^\alpha
\log\frac{2e^{1/\alpha}}{1-|\phi_t(z)|^2}
|f(\phi_t(z))|\notag\\
&\leq |w_t(z)|\,
\frac{\|f\|_{H_{\alpha,\log}^{\infty}}}
{(1-|\phi_t(z)|^2)^\alpha
\log\!\frac{2e^{1/\alpha}}{1-|\phi_t(z)|^2}}.
\label{eq:Hb-T-upper}
\end{align}
Noting that
\[
|w_t(z)|=\frac{(1-t)^b}{|1-(1-t)z|}
\leq\frac{(1-t)^b}{1-(1-t)|z|}=w_t(|z|),
\]
\[
|\phi_t(z)|=\frac{t}{|1-(1-t)z|}
\leq\frac{t}{1-(1-t)|z|}=\phi_t(|z|),
\]
and the monotonicity of the function
$$
g(x)=x^{\alpha}\log\frac{2e^{1/\alpha}}{x},~~ 0< x<2,
$$
 we obtain
\[
|T_tf(z)|\leq w_t(|z|)
\frac1{(1-\phi_t(|z|)^2)^\alpha
\log\!\frac{2e^{1/\alpha}}{1-\phi_t(|z|)^2}}
\|f\|_{H_{\alpha,\log}^{\infty}}.
\]
Therefore, we can obtain
\begin{align}
\|\mathcal H_bf\|_{H_\alpha^{\infty}}
&=\sup_{z\in\mathbb D}(1-|z|^2)^\alpha
\left|\int_0^1T_tf(z)\,dt\right|\notag\\
&\leq\sup_{0\leq r<1}(1-r^2)^\alpha
\int_0^1\frac{(1-t)^b}{1-(1-t)r}
\frac{dt}
{\left(1-\left(\frac{t}{1-(1-t)r}\right)^2\right)^\alpha
\log\!\frac{2e^{1/\alpha}}
{1-\left(\frac{t}{1-(1-t)r}\right)^2}}
\|f\|_{H_{\alpha,\log}^{\infty}}\notag\\
&=\sup_{0\leq r<1}\int_0^1
\frac{(1-t)^b(1+r)^\alpha\bigl((t-1)r+1\bigr)^{2\alpha-1}}
{(1-t)^\alpha\bigl((t-1)r+1+t\bigr)^\alpha
\log\!\frac{2e^{1/\alpha}}
{1-\left(\frac{t}{1-(1-t)r}\right)^2}}\,dt\;\|f\|_{H_{\alpha,\log}^{\infty}}.
\label{eq:Hb-upper}
\end{align}
Then we obtain
\begin{align}
&\|\mathcal H_b\|_{H_{\alpha,\log}^{\infty}\to H_\alpha^{\infty}}\notag\\
&\leq\sup_{0\leq r<1}\int_0^1
\frac{(1-t)^b(1+r)^\alpha\bigl((t-1)r+1\bigr)^{2\alpha-1}}
{(1-t)^\alpha\bigl((t-1)r+1+t\bigr)^\alpha
\log\!\frac{2e^{1/\alpha}}
{1-\left(\frac{t}{1-(1-t)r}\right)^2}}\,dt<\infty.
\label{eq:Hb-upper-norm}
\end{align}
The supremum is finite. Indeed, for $0\leq r<1$,
$(1-t)^b\leq(1-rt)^b$ and $f_\alpha(t)>0$ give
\[
\begin{aligned}
(1-r^2)^\alpha\mathcal H_bf_\alpha(r)
&\leq(1-r^2)^\alpha\int_0^1
\frac{f_\alpha(t)}{1-rt}\,dt\\
&\leq\|\mathcal H\|_{H_\alpha^\infty\to H_\alpha^\infty}
\|f_\alpha\|_{H_\alpha^\infty}<\infty.
\end{aligned}
\]
Here $\mathcal H$ is the classical Hilbert matrix operator, which
is bounded on $H_\alpha^\infty$ for $0<\alpha<1$,
and $\|f_\alpha\|_{H_\alpha^\infty}\leq
\|f_\alpha\|_{H_{\alpha,\log}^\infty}=1$.
Combining \eqref{eq:Hb-lower} and \eqref{eq:Hb-upper-norm}
proves the asserted norm equality.

Taking $r=0$ in the preceding supremum gives the stated lower bound.
 We finish the proof of
the theorem.
\end{proof}

Next, we recall the following geometric estimate, which will be used
in the proof of the next lemma. This estimate is given in
\cite{LindstromMiihkinenWikman2019}.

\begin{lemma}\label{lem:Hb-geometric}
Let $\frac12<\alpha<1$ and $0<t<1$. Then
\[
\begin{aligned}
&\sup_{z\in\mathbb D}|1-(1-t)z|^{2\alpha-1}
 \left(\frac{1-|z|^2}{|1-(1-t)z|^2-t^2}\right)^\alpha\\
&\qquad=\begin{cases}
\displaystyle\frac{t^{\alpha-1}}{(1-t)^\alpha},
&\begin{gathered}
\frac12<\alpha\leq\frac23,\ \text{or}\
\frac23<\alpha<1,\ t\geq\frac{3\alpha-2}{4\alpha-2},
\end{gathered}\\[4pt]
\displaystyle(1-x_0)^{2\alpha-1}
\left(\frac{1-\left|\frac{x_0}{1-t}\right|^2}
{(1-x_0)^2-t^2}\right)^\alpha,
&\begin{gathered}
\frac23<\alpha<1,\\
0<t<\frac{3\alpha-2}{4\alpha-2},
\end{gathered}
\end{cases}
\end{aligned}
\]
where
\[
x_0=\frac{\alpha+2\alpha t-t-
\sqrt{4\alpha^2t-2\alpha t+\alpha^2-2\alpha+1}}{2\alpha-1}.
\]
\end{lemma}

The preceding geometric estimate concerns only the symbol
$\phi_t(z)=t/[1-(1-t)z]$. Therefore, it is unchanged for the
operator $\mathcal H_b$; the additional factor $(1-t)^b$ appears
only in its weight $w_t$.

\begin{lemma}\label{lem:Hb-weighted-bound}
Let $b>0$, $0<\alpha<1$, and $0<t<1$. Then
\[
\begin{aligned}
&\|T_t\|_{H_{\alpha,\log}^{\infty}\to H_\alpha^\infty}\\[-2pt]
&\quad\leq\begin{cases}
\displaystyle
\frac{t^{\alpha-1}(1-t)^{b-\alpha}}
{\log\!\frac{(2-t)^2e^{1/\alpha}}{2(1-t)}},
&\begin{gathered}
0<\alpha\leq\frac23,\ \text{or}\
\frac23<\alpha<1,\ t\geq\frac{3\alpha-2}{4\alpha-2},
\end{gathered}\\[4pt]
\displaystyle
\frac{(1-t)^b(1-x_0)^{2\alpha-1}}
{\log\!\frac{(2-t)^2e^{1/\alpha}}{2(1-t)}}
\left(\frac{1-\left|\frac{x_0}{1-t}\right|^2}
{(1-x_0)^2-t^2}\right)^\alpha,
&\begin{gathered}
\frac23<\alpha<1,\\
0<t<\frac{3\alpha-2}{4\alpha-2},
\end{gathered}
\end{cases}
\end{aligned}
\]
where $x_0$ is as in Lemma~\ref{lem:Hb-geometric}.
\end{lemma}

\begin{proof}
Assume first that $0<\alpha\leq\frac12$ and
$f\in H_{\alpha,\log}^{\infty}$. Since
\[
\phi_t'(z)=\frac{t(1-t)}{[1-(1-t)z]^2},
\]
we have
\[
|\phi_t'(z)|\leq\frac{1-t}{t}.
\]
Moreover, the Schwarz--Pick lemma gives
\[
(1-|z|^2)|\phi_t'(z)|\leq1-|\phi_t(z)|^2.
\]
Because $|1-(1-t)z|\leq2-t$, we also have
\[
|\phi_t(z)|\geq\frac{t}{2-t},
\]
and hence
\[
\begin{aligned}
\frac1{\log\!\frac{2e^{1/\alpha}}
{1-|\phi_t(z)|^2}}
&\leq\frac1{\log\!\frac{2e^{1/\alpha}}
{1-\left(\frac{t}{2-t}\right)^2}}\\
&=\frac1{\log\!\frac{(2-t)^2e^{1/\alpha}}{2(1-t)}}.
\end{aligned}
\]
Using $|w_t(z)|=(1-t)^b|\phi_t'(z)|^{1/2}/\sqrt{t(1-t)}$,
we obtain
\begin{align}
&\|T_tf\|_{H_\alpha^\infty}
=\sup_{z\in\mathbb D}(1-|z|^2)^\alpha
|w_t(z)f(\phi_t(z))|\notag\\
&\leq\frac{(1-t)^b}{\sqrt{t(1-t)}}
\sup_{z\in\mathbb D}|\phi_t'(z)|^{1/2-\alpha}
\left[(1-|z|^2)|\phi_t'(z)|\right]^\alpha
|f(\phi_t(z))|\notag\\
&\leq\frac{(1-t)^b}{\sqrt{t(1-t)}}
\left(\frac{1-t}{t}\right)^{1/2-\alpha}
\sup_{z\in\mathbb D}
(1-|\phi_t(z)|^2)^\alpha|f(\phi_t(z))|\notag\\
&\leq\frac{t^{\alpha-1}(1-t)^{b-\alpha}}
{\log\!\frac{(2-t)^2e^{1/\alpha}}{2(1-t)}}
\|f\|_{H_{\alpha,\log}^\infty}.
\label{eq:Hb-lemma-small-alpha}
\end{align}
This gives the assertion for $0<\alpha\leq\frac12$.

Let $\frac12<\alpha<1$. As in \eqref{eq:Hb-T-upper},
\[
\begin{aligned}
&\|T_tf\|_{H_\alpha^\infty}\\
&\leq\frac{(1-t)^b}{\log\!\frac{(2-t)^2e^{1/\alpha}}{2(1-t)}}
\sup_{z\in\mathbb D}
\frac{(1-|z|^2)^\alpha}
{|1-(1-t)z|(1-|\phi_t(z)|^2)^\alpha}
\|f\|_{H_{\alpha,\log}^\infty}\\
&=\frac{(1-t)^b}{\log\!\frac{(2-t)^2e^{1/\alpha}}{2(1-t)}}
\sup_{z\in\mathbb D}|1-(1-t)z|^{2\alpha-1}
\left(\frac{1-|z|^2}{|1-(1-t)z|^2-t^2}\right)^\alpha
\|f\|_{H_{\alpha,\log}^\infty}.
\end{aligned}
\]
Applying Lemma~\ref{lem:Hb-geometric} gives the asserted
piecewise bound for $\frac12<\alpha<1$.
\end{proof}

\begin{theorem}\label{thm:Hb-integrated-upper}
Let $b>0$ and $0<\alpha<1$. If $0<\alpha\leq\frac23$,
then
\begin{equation}
\label{eq:Hb-simple-bound}
\|\mathcal H_b\|_{H_{\alpha,\log}^\infty\to H_\alpha^\infty}
\leq\int_0^1
\frac{t^{\alpha-1}(1-t)^{b-\alpha}}
{\log\!\frac{(2-t)^2e^{1/\alpha}}{2(1-t)}}\,dt.
\end{equation}
For $\frac23<\alpha<1$, then
\begin{equation}\label{eq:Hb-piecewise-bound}
\resizebox{0.91\linewidth}{!}{$\displaystyle
\|\mathcal H_b\|_{H_{\alpha,\log}^\infty\to H_\alpha^\infty}
\leq\int_0^{\frac{3\alpha-2}{4\alpha-2}}
\frac{(1-t)^bG(x_0)}{\log\!\frac{(2-t)^2e^{1/\alpha}}{2(1-t)}}\,dt
+\int_{\frac{3\alpha-2}{4\alpha-2}}^1
\frac{t^{\alpha-1}(1-t)^{b-\alpha}}{\log\!\frac{(2-t)^2e^{1/\alpha}}{2(1-t)}}\,dt.
$}
\end{equation}
where
\[
G(x)=(1-x)^{2\alpha-1}
\left(\frac{1-\left|\frac{x}{1-t}\right|^2}
{(1-x)^2-t^2}\right)^\alpha,
\]
and
\[
x_0=\frac{\alpha+2\alpha t-t-
\sqrt{4\alpha^2t-2\alpha t+\alpha^2-2\alpha+1}}{2\alpha-1}.
\]
All the integrals on the right-hand sides are finite.
\end{theorem}

\begin{proof}
Let $b>0$ and $0<\alpha<1$. By the weighted composition
representation \eqref{eq:Hb-weighted-composition}, we have
\begin{align}
\|\mathcal H_bf\|_{H_\alpha^\infty}
&=\sup_{z\in\mathbb D}(1-|z|^2)^\alpha
\left|\int_0^1T_tf(z)\,dt\right|\notag\\
&\leq\sup_{z\in\mathbb D}\int_0^1
(1-|z|^2)^\alpha|T_tf(z)|\,dt\notag\\
&\leq\int_0^1\sup_{z\in\mathbb D}
(1-|z|^2)^\alpha|T_tf(z)|\,dt\notag\\
&=\int_0^1\|T_tf\|_{H_\alpha^\infty}\,dt\notag\\
&\leq\int_0^1
\|T_t\|_{H_{\alpha,\log}^\infty\to H_\alpha^\infty}
\,dt\;\|f\|_{H_{\alpha,\log}^\infty}.
\label{eq:Hb-minkowski}
\end{align}
Combining \eqref{eq:Hb-minkowski} with
Lemma~\ref{lem:Hb-weighted-bound} gives
\eqref{eq:Hb-simple-bound} for $0<\alpha\leq\frac23$
and \eqref{eq:Hb-piecewise-bound} for $\frac23<\alpha<1$.

To see that the bounds are finite, observe first that the
logarithmic denominator is bounded below by
$\log(2e^{1/\alpha})>0$. For $0<\alpha\leq\frac23$,
the integrand in \eqref{eq:Hb-simple-bound} is bounded by
\[
\frac{t^{\alpha-1}(1-t)^{b-\alpha}}
{\log(2e^{1/\alpha})},
\]
which is integrable because $\alpha>0$ and
$b-\alpha>-1$.
For $\frac23<\alpha<1$, the second integral is handled in
the same way. For the first integral, the Schwarz--Pick estimate
used above gives
\[
|1-(1-t)z|^2-t^2
\geq t(1-t)(1-|z|^2),\qquad z\in\mathbb D.
\]
Since $2\alpha-1>0$ and
$|1-(1-t)z|\leq 2-t$, Lemma~\ref{lem:Hb-geometric}
therefore implies
\[
G(x_0)\leq\frac{(2-t)^{2\alpha-1}}{[t(1-t)]^\alpha}
\qquad\left(0<t<\frac{3\alpha-2}{4\alpha-2}\right).
\]
Consequently, the first integrand in
\eqref{eq:Hb-piecewise-bound} is at most
\[
\frac{(1-t)^{b-\alpha}(2-t)^{2\alpha-1}}
{t^\alpha\log(2e^{1/\alpha})}.
\]
Since $\alpha<1$, the factor $t^{-\alpha}$ is integrable
near $t=0$. Moreover, for $\frac23<\alpha<1$,
\[
0<\frac{3\alpha-2}{4\alpha-2}<1,
\]
so the remaining factors are bounded on the interval away from
$t=0$. Thus the first integral is finite. This finishes the
proof of the theorem.
\end{proof}

\begin{remark}
Setting $b=0$, the above two theorems reduce to Theorems 3.1 and 3.2 in \cite{HuYe2025}, respectively.
\end{remark}

\section[Norm estimates on Bloch-type spaces]{Norm estimates for the generalized Hilbert matrix on $\mathcal B^\alpha$}

In this section, we study norm estimates for the generalized Hilbert
matrix operator $\mathcal H_b$ acting on $\mathcal B^\alpha$. We first
consider $1<\alpha<2$, as in the classical case, and then indicate the
larger range allowed by the factor $(1-t)^b$ when $b>0$. For convenience,
we write
\[
\|f\|_{\alpha,*}=\sup_{z\in\mathbb D}(1-|z|^2)^\alpha|f'(z)|.
\]
By differentiating the integral representation
\eqref{eq:Hb-integral}, we obtain
\begin{equation}\label{eq:Hb-Bloch-derivative}
(\mathcal H_bf)'(z)
=(b+1)\int_0^1\frac{t(1-t)^b}{(1-tz)^{b+2}}f(t)\,dt.
\end{equation}
The differentiation is justified locally uniformly in $z$ whenever
$\int_0^1(1-t)^b|f(t)|\,dt<\infty$. For $z\in\mathbb D$, consider the
same map as in Section~\ref{sec:Hb-representation},
\[
\phi_t(z)=\frac{t}{(t-1)z+1}.
\]
For real $0\leq z<1$, the change of variables in
\eqref{eq:Hb-Bloch-derivative} gives
\begin{equation}\label{eq:Hb-Bloch-composition-derivative}
(\mathcal H_bf)'(z)
=(b+1)\int_0^1
\frac{t(1-t)^b}{[(t-1)z+1](1-z)}
f\bigl(\phi_t(z)\bigr)\,dt.
\end{equation}
Both sides are analytic functions of $z\in\mathbb D$; hence
\eqref{eq:Hb-Bloch-composition-derivative} holds throughout $\mathbb D$.
The factor $(b+1)(1-t)^b$ is essential here; for $b=0$, the
representation reduces to the derivative formula for the classical
Hilbert matrix operator.

We recall the following radial growth estimate from Hu and Ye
\cite[Lemma~5.1]{HuYe2025}.

\begin{lemma}\label{lem:Hb-Bloch-growth}
Let $\alpha>0$ and $f\in\mathcal B^\alpha$. Then
\[
|f(z)|\leq
\begin{cases}
\displaystyle |f(0)|+\frac{(1-|z|)^{1-\alpha}-1}{\alpha-1}
\|f\|_{\alpha,*},&\alpha\ne1,\\[7pt]
\displaystyle |f(0)|+\log\frac1{1-|z|}\,\|f\|_{\alpha,*},
&\alpha=1.
\end{cases}
\]
\end{lemma}

\begin{theorem}\label{thm:Hb-Bloch-lower}
Let $b>0$ and $1<\alpha<b+2$. Then
\begin{equation}\label{eq:Hb-Bloch-lower}
\|\mathcal H_b\|_{\mathcal B^\alpha\to\mathcal B^\alpha}
\geq\frac{\displaystyle\int_0^1(1-t)^b(1-t^2)^{1-\alpha}\,dt-\frac1{b+1}}{2(\alpha-1)}
+\frac{b+1}{\alpha-1}B(b+2-\alpha,\alpha).
\end{equation}
Moreover, $\mathcal H_b$ is not a bounded operator on
$\mathcal B^\alpha$ when $0<\alpha\leq1$ or
$\alpha\geq b+2$.
\end{theorem}

\begin{proof}
Let $\alpha>1$ and define the same test function as in the
classical argument,
\[
h_\alpha(z)=\frac{(1-z^2)^{1-\alpha}-1}{2(\alpha-1)}.
\]
Since $1-z^2$ lies in the right half-plane for $z\in\mathbb D$,
the analytic branch in this formula is well defined. Direct
differentiation gives
\[
h_\alpha(0)=0,\qquad
h_\alpha'(z)=\frac{z}{(1-z^2)^\alpha}.
\]
On the one hand,
\[
\|h_\alpha\|_{\mathcal B^\alpha}
=\sup_{z\in\mathbb D}\frac{|z|(1-|z|^2)^\alpha}
{|1-z^2|^\alpha}\leq1.
\]
On the other hand, for real $0<r<1$,
\[
(1-r^2)^\alpha|h_\alpha'(r)|=r,
\]
and letting $r$ tend to $1$ from below gives
\[
\|h_\alpha\|_{\mathcal B^\alpha}=1.
\]
For $1<\alpha<b+2$, the defining integrals for
$\mathcal H_bh_\alpha$ converge. By
\eqref{eq:Hb-integral},
\[
\begin{aligned}
(\mathcal H_bh_\alpha)(0)
&=\int_0^1(1-t)^bh_\alpha(t)\,dt\\
&=\frac1{2(\alpha-1)}
\int_0^1(1-t)^b(1-t^2)^{1-\alpha}\,dt
-\frac1{2(\alpha-1)(b+1)}.
\end{aligned}
\]
In addition, by \eqref{eq:Hb-Bloch-composition-derivative},
\[
(\mathcal H_bh_\alpha)'(z)
=\frac{b+1}{2(\alpha-1)}\int_0^1
\frac{t(1-t)^b}{[(t-1)z+1](1-z)}
\left[\left(1-\phi_t(z)^2\right)^{1-\alpha}-1\right]dt.
\]
In particular, for real $0\leq r<1$, the integrand is
nonnegative and, using \eqref{eq:Hb-Bloch-derivative},
\[
\begin{aligned}
&(1-r^2)^\alpha(\mathcal H_bh_\alpha)'(r)\\
&\quad=\frac{b+1}{2(\alpha-1)}(1-r^2)^\alpha
\int_0^1\frac{t(1-t)^b(1-t^2)^{1-\alpha}}
{(1-rt)^{b+2}}\,dt\\
&\qquad-\frac{b+1}{2(\alpha-1)}(1-r^2)^\alpha
\int_0^1\frac{t(1-t)^b}{(1-rt)^{b+2}}\,dt.
\end{aligned}
\]
We first determine the limit of the second integral term. The
substitution $x=(1-t)/(1-rt)$ gives
\begin{equation}\label{eq:Hb-Bloch-elementary-integral}
\int_0^1\frac{(1-t)^b}{(1-rt)^{b+2}}\,dt
=\frac{1}{(b+1)(1-r)}.
\end{equation}
Consequently,
\[
0\leq(b+1)(1-r^2)^\alpha
\int_0^1\frac{t(1-t)^b}{(1-rt)^{b+2}}\,dt
\leq\frac{(1-r^2)^\alpha}{1-r}
\longrightarrow0
\]
as $r$ tends to $1$ from below, because $\alpha>1$.

For the first integral term, put
$u=(1-t)/(1-r)$. A direct calculation gives
\begin{equation}\label{eq:Hb-Bloch-boundary-substitution}
\begin{aligned}
&(1-r^2)^\alpha\int_0^1
\frac{t(1-t)^b(1-t^2)^{1-\alpha}}
{(1-rt)^{b+2}}\,dt\\
&\quad=(1+r)^\alpha
\int_0^{1/(1-r)}
\frac{[1-(1-r)u]u^{b+1-\alpha}
[2-(1-r)u]^{1-\alpha}}
{(1+ru)^{b+2}}\,du.
\end{aligned}
\end{equation}
For fixed $u>0$, the integrand on the right tends to
$2u^{b+1-\alpha}/(1+u)^{b+2}$. If $r\geq1/2$, its
absolute value, extended by zero outside the integration interval,
is bounded by
\[
2^\alpha\frac{u^{b+1-\alpha}}{(1+u/2)^{b+2}},
\]
which is integrable on $(0,\infty)$ because
$1<\alpha<b+2$. Hence dominated convergence in
\eqref{eq:Hb-Bloch-boundary-substitution} yields
\[
\begin{aligned}
&\lim_{r\to1^-}(1-r^2)^\alpha
\int_0^1\frac{t(1-t)^b(1-t^2)^{1-\alpha}}
{(1-rt)^{b+2}}\,dt\\
&\qquad=2\int_0^\infty
\frac{u^{b+1-\alpha}}{(1+u)^{b+2}}\,du
=2B(b+2-\alpha,\alpha).
\end{aligned}
\]
Therefore,
\begin{equation}\label{eq:Hb-Bloch-boundary-derivative}
\lim_{r\to1^-}(1-r^2)^\alpha
(\mathcal H_bh_\alpha)'(r)
=\frac{b+1}{\alpha-1}B(b+2-\alpha,\alpha).
\end{equation}
Since $\|h_\alpha\|_{\mathcal B^\alpha}=1$ and
$(\mathcal H_bh_\alpha)(0)>0$, we obtain
\[
\begin{aligned}
\|\mathcal H_b\|_{\mathcal B^\alpha\to\mathcal B^\alpha}
&\geq\|\mathcal H_bh_\alpha\|_{\mathcal B^\alpha}\\
&\geq(\mathcal H_bh_\alpha)(0)
+\lim_{r\to1^-}(1-r^2)^\alpha
(\mathcal H_bh_\alpha)'(r),
\end{aligned}
\]
and \eqref{eq:Hb-Bloch-lower} follows.

We next consider the remaining parameters. If $0<\alpha<1$,
take $f\equiv1$, which has $\|f\|_{\mathcal B^\alpha}=1$.
For real $r=1-\delta$, with $0<\delta<1/4$, the interval
$1-2\delta\leq t\leq1-\delta$ is contained in $(0,1)$, and
\[
\begin{aligned}
&(1-r^2)^\alpha(\mathcal H_b1)'(r)\\
&\quad=(b+1)(1-r^2)^\alpha
\int_0^1\frac{t(1-t)^b}{(1-rt)^{b+2}}\,dt\\
&\quad\geq \frac{b+1}{2\cdot3^{b+2}}\,
\delta^{\alpha-1}\longrightarrow\infty.
\end{aligned}
\]
Thus $\mathcal H_b$ is not bounded on $\mathcal B^\alpha$
for $0<\alpha<1$.

For $\alpha=1$, let $h_1(z)=\log\frac1{1-z}$.
Since $h_1(0)=0$ and
\[
(1-|z|^2)|h_1'(z)|
=\frac{1-|z|^2}{|1-z|}\leq1+|z|\leq2,
\]
we have $h_1\in\mathcal B^1$. On the same interval
$1-2\delta\leq t\leq1-\delta$, the derivative representation
gives
\[
\begin{aligned}
&(1-r^2)(\mathcal H_bh_1)'(r)\\
&\quad=(b+1)(1-r^2)
\int_0^1\frac{t(1-t)^b}{(1-rt)^{b+2}}
\log\frac1{1-t}\,dt\\
&\quad\geq\frac{b+1}{2\cdot3^{b+2}}
\log\frac1{2\delta}\longrightarrow\infty.
\end{aligned}
\]
Hence $\mathcal H_b$ is not bounded on $\mathcal B^1$.

Finally, if $\alpha\geq b+2$, use the same function
$h_\alpha$ as above. We still have
$\|h_\alpha\|_{\mathcal B^\alpha}=1$, but
\[
(\mathcal H_bh_\alpha)(0)
=\frac1{2(\alpha-1)}
\int_0^1(1-t)^b\bigl[(1-t^2)^{1-\alpha}-1\bigr]dt
=\infty,
\]
because the first integrand is comparable to
$(1-t)^{b+1-\alpha}$ near $t=1$, and
$b+1-\alpha\leq-1$. Therefore $\mathcal H_b$ cannot be
a bounded operator on $\mathcal B^\alpha$ in this range.
\end{proof}

Next, we consider the upper bound for the generalized Hilbert
matrix operator acting on $\mathcal B^\alpha$.

\begin{theorem}\label{thm:Hb-Bloch-upper}
Let $b>0$ and $1<\alpha<b+2$. Then $\mathcal H_b$ is bounded on
$\mathcal B^\alpha$, and its norm satisfies
\begin{equation}\label{eq:Hb-Bloch-upper}
\|\mathcal H_b\|_{\mathcal B^\alpha\to\mathcal B^\alpha}
\leq\max\!\left\{\frac{2b+3}{b+1},\;
\frac{1}{(b+1)(b+2-\alpha)}
+\frac{2^\alpha(b+1)}{\alpha-1}
B(b+2-\alpha,\min\{\alpha,2\})\right\}.
\end{equation}
Consequently, for $b>0$ and $\alpha>0$, the operator
$\mathcal H_b:\mathcal B^\alpha\to\mathcal B^\alpha$ is bounded
if and only if $1<\alpha<b+2$.
\end{theorem}

\begin{proof}
Let $1<\alpha<b+2$ and $f\in\mathcal B^\alpha$.
According to Lemma~\ref{lem:Hb-Bloch-growth},
\[
|f(t)|\leq|f(0)|+
\frac{(1-t)^{1-\alpha}-1}{\alpha-1}\|f\|_{\alpha,*},
\qquad0\leq t<1.
\]
By \eqref{eq:Hb-integral}, we obtain
\[
\begin{aligned}
|\mathcal H_bf(0)|
&\leq\int_0^1(1-t)^b|f(t)|\,dt\\
&\leq\frac1{b+1}|f(0)|
+\frac1{\alpha-1}
\left(\frac1{b+2-\alpha}-\frac1{b+1}\right)
\|f\|_{\alpha,*}\\
&=\frac1{b+1}|f(0)|
+\frac1{(b+1)(b+2-\alpha)}\|f\|_{\alpha,*}.
\end{aligned}
\]
For the derivative seminorm, \eqref{eq:Hb-Bloch-derivative}
and $|1-tz|\geq1-t|z|$ imply, with $r=|z|$,
\[
\begin{aligned}
&(1-|z|^2)^\alpha|(\mathcal H_bf)'(z)|\\
&\quad\leq(b+1)(1-r^2)^\alpha
\int_0^1\frac{t(1-t)^b}{(1-rt)^{b+2}}
\left(|f(0)|+
\frac{(1-t)^{1-\alpha}-1}{\alpha-1}
\|f\|_{\alpha,*}\right)dt.
\end{aligned}
\]
For the term containing $|f(0)|$, we use
\eqref{eq:Hb-Bloch-elementary-integral} and $t\leq1$:
\[
\begin{aligned}
&(b+1)(1-r^2)^\alpha
\int_0^1\frac{t(1-t)^b}{(1-rt)^{b+2}}\,dt\\
&\qquad\leq\frac{(1-r^2)^\alpha}{1-r}
=(1+r)(1-r^2)^{\alpha-1}\leq1+r<2.
\end{aligned}
\]
For the term containing $\|f\|_{\alpha,*}$, discard the
negative term in the growth bound and use the substitution
\[
x=\frac{1-t}{1-rt},\qquad
 t=\frac{1-x}{1-rx}.
\]
Then
\begin{equation}\label{eq:Hb-Bloch-upper-substitution}
\begin{aligned}
&(1-r^2)^\alpha
\int_0^1\frac{t(1-t)^{b+1-\alpha}}
{(1-rt)^{b+2}}\,dt\\
&\quad=(1+r)^\alpha\int_0^1
(1-x)x^{b+1-\alpha}(1-rx)^{\alpha-2}\,dx.
\end{aligned}
\end{equation}
If $1<\alpha\leq2$, then $\alpha-2\leq0$ and
$1-rx\geq1-x$, so
\[
(1-x)(1-rx)^{\alpha-2}\leq(1-x)^{\alpha-1}.
\]
Consequently, \eqref{eq:Hb-Bloch-upper-substitution} yields
\[
\begin{aligned}
&(1-r^2)^\alpha
\int_0^1\frac{t(1-t)^{b+1-\alpha}}
{(1-rt)^{b+2}}\,dt\\
&\quad\leq(1+r)^\alpha
\int_0^1x^{b+1-\alpha}(1-x)^{\alpha-1}\,dx\\
&\quad\leq2^\alpha B(b+2-\alpha,\alpha).
\end{aligned}
\]
If $2<\alpha<b+2$, then $0<\alpha-2$ and
$1-rx\leq1$, whence
\[
(1-x)(1-rx)^{\alpha-2}\leq1-x.
\]
The same calculation gives
\[
\begin{aligned}
&(1-r^2)^\alpha
\int_0^1\frac{t(1-t)^{b+1-\alpha}}
{(1-rt)^{b+2}}\,dt\\
&\quad\leq2^\alpha\int_0^1
x^{b+1-\alpha}(1-x)\,dx
=2^\alpha B(b+2-\alpha,2).
\end{aligned}
\]
Combining the estimates, we obtain
\[
\resizebox{\linewidth}{!}{$\displaystyle
\|\mathcal H_bf\|_{\mathcal B^\alpha}
\leq \left(2+\frac1{b+1}\right)|f(0)|
+\left[\frac1{(b+1)(b+2-\alpha)}
+\frac{2^\alpha(b+1)}{\alpha-1}
B\bigl(b+2-\alpha,\min\{\alpha,2\}\bigr)\right]\|f\|_{\alpha,*}.
$}
\]
Since $\|f\|_{\mathcal B^\alpha}=|f(0)|+\|f\|_{\alpha,*}$,
we conclude that
\eqref{eq:Hb-Bloch-upper} holds, proving boundedness in
$1<\alpha<b+2$. The assertion on the full boundedness range
follows from Theorem~\ref{thm:Hb-Bloch-lower}.
\end{proof}
\begin{remark}
Setting \(b=0\) in Theorems 4.2 and 4.3, we obtain Theorems 5.1 and 5.2 in \cite{HuYe2025}, respectively.
\end{remark}

\section[Norm estimates from $H^\infty$ to a Zygmund-type space]{Norm estimates for $\mathcal H_b$ from $H^\infty$ to $\Lambda_1^{1,*}$}

Li, Pan and Zhou \cite{LiPanZhou2026} recently determined the norm of
$\mathcal H_b$ from $H^\infty$ into the classical Bloch space:
\[
\|\mathcal H_b\|_{H^\infty\to\mathcal B}
=\frac{1}{b+1}+2,\qquad b>0.
\]
The Bloch norm involves the first derivative, whereas the endpoint
mean Lipschitz space considered by Guo and Tang \cite{GuoTang2024}
is described in terms of integral means of the second derivative.
We now consider the corresponding mapping of $\mathcal H_b$ into this  endpoint space
 and establish two-sided norm estimates.

Following \cite{GuoTang2024}, we write $\Lambda_1^{1,*}=\mathcal Z_1$, equipped with the norm
\[
\|f\|_{\Lambda_1^{1,*}}
=|f(0)|+|f'(0)|+
\sup_{0\leq r<1}(1-r^2)M_1(r,f''),
\]
where $M_1$ is defined in Section~2. From
\eqref{eq:Hb-integral}, differentiation under the integral sign gives
\begin{equation}\label{eq:Hb-Zygmund-second-derivative}
(\mathcal H_b f)''(z)
=(b+1)(b+2)\int_0^1
\frac{t^2(1-t)^b f(t)}{(1-tz)^{b+3}}\,dt.
\end{equation}
The integral and its derivatives converge locally uniformly on
$\mathbb D$ for $f\in H^\infty$.

We also need the following estimates, which are the one-dimensional
case of \cite[Theorem~1.3]{Liu2015}.
\begin{lemma}\label{lem:Liu-Forelli-Rudin}
For $z\in\mathbb D$ and $c\in\mathbb R$, define
\[
I_c(z):=\frac1{2\pi}\int_0^{2\pi}
\frac{d\theta}{|1-ze^{-i\theta}|^{1+c}}.
\]
Then the following statements hold:

(1)  If $c<0$, then 
$$ 1\leq I_c(z)\leq
\frac{\Gamma(-c)}{\Gamma^2\!\left(\frac{1-c}{2}\right)}.
$$
(2)  If $c=0$, then 
\[
1\leq(1-|z|^2)^c I_c(z)\leq
\frac{\Gamma(c)}{\Gamma^2\!\left(\frac{1+c}{2}\right)}.
\]
(1)  If $c>$, then 
\[
\frac1\pi\leq |z|^2\left(\log\frac1{1-|z|^2}\right)^{-1}I_0(z)\leq1.
\]
Moreover, all these inequalities are sharp.
\end{lemma}

\begin{theorem}\label{thm:Hb-Zygmund-bounds}
Let $b>0$. Then the generalized Hilbert operator
$\mathcal H_b:H^\infty\to\Lambda_1^{1,*}$ is bounded, and
\begin{equation}\label{eq:Hb-Zygmund-two-sided}
\frac{2b+3}{(b+1)(b+2)}+\frac2\pi
\leq\|\mathcal H_b\|_{H^\infty\to\Lambda_1^{1,*}}
\leq\frac{2b+3}{(b+1)(b+2)}
+\frac{(b+1)(b+2)\Gamma(b+2)}{2\Gamma^2\!\left(\frac{b+3}{2}\right)}.
\end{equation}
\end{theorem}

\begin{proof}
We first prove the lower bound. Let $f(z)=1$. Then
$\|f\|_{H^\infty}=1$, and the integral representation gives
\begin{equation}\label{eq:Hb-Zygmund-origin-values}
\mathcal H_b f(0)=\int_0^1(1-t)^b\,dt=\frac{1}{b+1},
\qquad
(\mathcal H_b f)'(0)
=(b+1)\int_0^1t(1-t)^b\,dt=\frac{1}{b+2}.
\end{equation}
For any analytic function $g$, the Fej\'er--Riesz inequality
(see \cite{Duren1970}) implies
\[
\int_{-1}^1 |g(rx)|\,dx\leq\pi M_1(r,g),
\qquad 0<r<1.
\]
By \eqref{eq:Hb-Zygmund-second-derivative},
$(\mathcal H_b f)''(x)>0$ for real $0\leq x<1$.
Consequently,
\begin{align}
M_1\bigl(r,(\mathcal H_b f)''\bigr)
&\geq\frac1\pi\int_0^1(\mathcal H_b f)''(rx)\,dx\notag\\
&=\frac{(\mathcal H_b f)'(r)
-(\mathcal H_b f)'(0)}{\pi r}.
\label{eq:Hb-Zygmund-lower-radial}
\end{align}
To calculate its boundary contribution, observe that
\[
(\mathcal H_b f)'(r)
=(b+1)\int_0^1
\frac{t(1-t)^b}{(1-rt)^{b+2}}\,dt.
\]
Putting $t=1-(1-r)x$, we obtain
\begin{equation}\label{eq:Hb-Zygmund-boundary-computation}
\begin{aligned}
(1-r)(\mathcal H_b f)'(r)
&=(b+1)\int_0^{1/(1-r)}
\frac{[1-(1-r)x]x^b}{(1+rx)^{b+2}}\,dx,\\
\lim_{r\to1^-}(1-r)(\mathcal H_b f)'(r)
&=(b+1)\int_0^\infty\frac{x^b}{(1+x)^{b+2}}\,dx\\
&=(b+1)B(b+1,1)=1.
\end{aligned}
\end{equation}
Indeed, for $r\geq\tfrac12$ the integrand, extended by
zero outside $[0,1/(1-r)]$, is bounded by
$x^b/(1+x/2)^{b+2}$, which belongs to $L^1(0,\infty)$.
Thus the limit follows from dominated convergence. Since
$(\mathcal H_b f)'(0)$ is finite, \eqref{eq:Hb-Zygmund-lower-radial}
and \eqref{eq:Hb-Zygmund-boundary-computation} yield
\[
\begin{aligned}
&\sup_{0<r<1}(1-r^2)
M_1\bigl(r,(\mathcal H_b f)''\bigr)\\
&\qquad\geq \frac1\pi\lim_{r\to1^-}
\frac{1-r^2}{r}
\bigl((\mathcal H_b f)'(r)-(\mathcal H_b f)'(0)\bigr)
=\frac2\pi.
\end{aligned}
\]
Combining this with \eqref{eq:Hb-Zygmund-origin-values},
we obtain
\[
\|\mathcal H_b\|_{H^\infty\to\Lambda_1^{1,*}}
\geq\|\mathcal H_b f\|_{\Lambda_1^{1,*}}
\geq\frac1{b+1}+\frac1{b+2}+\frac2\pi.
\]

We next consider the upper bound. For $f\in H^\infty$,
\eqref{eq:Hb-integral} gives
\[
\left|\mathcal H_bf(0)\right|
\leq\frac1{b+1}\|f\|_{H^\infty},
\qquad
\left|(\mathcal H_bf)'(0)\right|
\leq\frac1{b+2}\|f\|_{H^\infty}.
\]
To estimate the second derivative, we apply part (2) of
Lemma~\ref{lem:Liu-Forelli-Rudin} with $c=b+2>0$ and $z=q$.
Since changing $e^{-i\theta}$ to $e^{i\theta}$ does not affect
the angular integral, we obtain, for $b>0$ and $0\leq q<1$,
\begin{equation}\label{eq:Hb-Zygmund-angular-bound}
\frac1{2\pi}\int_0^{2\pi}
\frac{d\theta}{|1-qe^{i\theta}|^{b+3}}
\leq
\frac{\Gamma(b+2)}
 {\Gamma^2\!\left(\frac{b+3}{2}\right)}
\frac1{(1-q^2)^{b+2}}.
\end{equation}

By \eqref{eq:Hb-Zygmund-second-derivative}, Minkowski's
integral inequality, and \eqref{eq:Hb-Zygmund-angular-bound},
\[
\begin{aligned}
&M_1\bigl(r,(\mathcal H_b f)''\bigr)\\
&\quad\leq(b+1)(b+2)\|f\|_{H^\infty}
\int_0^1t^2(1-t)^b
\left(\frac1{2\pi}\int_0^{2\pi}
\frac{d\theta}{|1-tre^{i\theta}|^{b+3}}\right)dt\\
&\quad\leq\frac{(b+1)(b+2)\Gamma(b+2)}
 {\Gamma^2\!\left(\frac{b+3}{2}\right)}
\|f\|_{H^\infty}
\int_0^1\frac{t^2(1-t)^b}{(1-t^2r^2)^{b+2}}\,dt.
\end{aligned}
\]
Since $1-t\leq1-t^2r^2$ for $0\leq r,t<1$ and $b\geq0$,
\[
\frac{(1-t)^b}{(1-t^2r^2)^b}\leq1.
\]
Also $t^2\leq t$ on $[0,1]$. Thus
\begin{equation}\label{eq:Hb-Zygmund-upper-radial}
\begin{aligned}
&(1-r^2)\int_0^1
\frac{t^2(1-t)^b}{(1-t^2r^2)^{b+2}}\,dt\\
&\quad\leq(1-r^2)\int_0^1
\frac{t}{(1-t^2r^2)^2}\,dt=\frac12.
\end{aligned}
\end{equation}
For $r=0$ the final equality is immediate; for $0<r<1$
it follows by direct integration. Consequently,
\[
\begin{aligned}
&\sup_{0\leq r<1}(1-r^2)
M_1\bigl(r,(\mathcal H_bf)''\bigr)\\
&\quad\leq
\frac{(b+1)(b+2)\Gamma(b+2)}
 {2\Gamma^2\!\left(\frac{b+3}{2}\right)}
\|f\|_{H^\infty}.
\end{aligned}
\]
Combining this estimate with those for the two terms at the
origin, we obtain
\[
\begin{aligned}
\|\mathcal H_bf\|_{\Lambda_1^{1,*}}
&\leq\left[\frac1{b+1}+\frac1{b+2}
+\frac{(b+1)(b+2)\Gamma(b+2)}
 {2\Gamma^2\!\left(\frac{b+3}{2}\right)}\right]
\|f\|_{H^\infty}.
\end{aligned}
\]
This proves the upper bound in \eqref{eq:Hb-Zygmund-two-sided}
and the boundedness of $\mathcal H_b:H^\infty\to\Lambda_1^{1,*}$.
\end{proof}

\begin{remark}
Setting \(b=0\) in Theorems 5.2, we obtain Theorems 1.2  in \cite{GuoTang2024}.
\end{remark}

\section[Norm estimates from $H_\alpha^\infty$ to $\mathcal B^{\alpha+1}$]{Norm estimates for $\mathcal H_b$ from $H_\alpha^\infty$ to $\mathcal B^{\alpha+1}$}

Hu and Ye \cite[Theorem~7.1]{HuYe2025} studied the Hilbert matrix
operator from the Korenblum space $H_\alpha^\infty$ into the
$(\alpha+1)$-Bloch space $\mathcal B^{\alpha+1}$. We consider here the corresponding mapping for the generalized Hilbert matrix operator
$\mathcal H_b$. The factor $(1-t)^b$ in its integral kernel changes
the admissible range of $\alpha$, while the estimates can be obtained
by a similar argument.

According to \eqref{eq:Hb-integral}, for $f\in H_\alpha^\infty$ and
$0<\alpha<b+1$ we have
\begin{equation}\label{eq:Hb-Korenblum-Bloch-derivative}
(\mathcal H_bf)'(z)
=(b+1)\int_0^1\frac{t(1-t)^b}{(1-tz)^{b+2}}f(t)\,dt,
\qquad z\in\mathbb D.
\end{equation}
Here and below the powers of $1-tz$ are defined by the analytic
branch that is positive on $z\in[0,1)$. We use the test function
\[
f_\alpha(z)=\frac{1}{(1-z^2)^\alpha}.
\]

\begin{theorem}\label{thm:Hb-Korenblum-Bloch}
Let $b>0$. For $0<\alpha<b+1$, the operator
$\mathcal H_b:H_\alpha^\infty\to\mathcal B^{\alpha+1}$
is bounded and
\begin{equation}\label{eq:Hb-Korenblum-Bloch-norm}
\resizebox{0.91\linewidth}{!}{$\displaystyle
\|\mathcal H_b\|_{H_\alpha^\infty\to\mathcal B^{\alpha+1}}
=\int_0^1\frac{(1-t)^b}{(1-t^2)^\alpha}\,dt
+(b+1)\sup_{0\leq r<1}(1-r^2)^{\alpha+1}
\int_0^1\frac{t(1-t)^b}{(1-rt)^{b+2}(1-t^2)^\alpha}\,dt.
$}
\end{equation}
Moreover, $f_\alpha(z)=(1-z^2)^{-\alpha}$ has norm one in
$H_\alpha^\infty$ and attains the norm of  $\mathcal H_b$.

\noindent For $0<\alpha\leq\frac23$, we have
\begin{equation}\label{eq:Hb-Korenblum-Bloch-sharp}
\|\mathcal H_b\|_{H_\alpha^\infty\to\mathcal B^{\alpha+1}}
=\int_0^1\frac{(1-t)^b}{(1-t^2)^\alpha}\,dt
 +2(b+1)B(b+1-\alpha,1+\alpha).
\end{equation}

\noindent For $\frac23<\alpha<\min\{2,b+1\}$, the following lower and
upper bounds hold:
\begin{equation}\label{eq:Hb-Korenblum-Bloch-bounds}
\resizebox{0.99\linewidth}{!}{$\displaystyle
\int_0^1\frac{(1-t)^b}{(1-t^2)^\alpha}\,dt
+2(b+1)B(b+1-\alpha,1+\alpha)
\leq\|\mathcal H_b\|_{H_\alpha^\infty\to\mathcal B^{\alpha+1}}
\leq\int_0^1\frac{(1-t)^b}{(1-t^2)^\alpha}\,dt
+2(b+1)B(b+1-\alpha,2-\alpha).
$}
\end{equation}
If $\alpha\geq b+1$, then $\mathcal H_b$ does not define a bounded
operator from $H_\alpha^\infty$ into $\mathcal B^{\alpha+1}$.
\end{theorem}

\begin{proof}
First, we consider the lower bound and the norm reduction. For
$0<\alpha<b+1$, the test function $f_\alpha$ belongs to
$H_\alpha^\infty$. Indeed, using $|1-z^2|\geq1-|z|^2$, we have
\[
\|f_\alpha\|_{H_\alpha^\infty}
=\sup_{z\in\mathbb D}\frac{(1-|z|^2)^\alpha}{|1-z^2|^\alpha}
\leq1.
\]
Equality holds on the real interval $[0,1)$, and hence
$\|f_\alpha\|_{H_\alpha^\infty}=1$. By
\eqref{eq:Hb-integral} and \eqref{eq:Hb-Korenblum-Bloch-derivative},
\[
\mathcal H_bf_\alpha(0)
=\int_0^1\frac{(1-t)^b}{(1-t^2)^\alpha}\,dt
\]
and, for $0\leq r<1$,
\[
(\mathcal H_bf_\alpha)'(r)
=(b+1)\int_0^1\frac{t(1-t)^b}{(1-rt)^{b+2}(1-t^2)^\alpha}\,dt>0.
\]
Consequently,
\begin{equation}\label{eq:Hb-Korenblum-Bloch-lower}
\resizebox{0.91\linewidth}{!}{$\displaystyle
\|\mathcal H_b\|_{H_\alpha^\infty\to\mathcal B^{\alpha+1}}
\geq\|\mathcal H_bf_\alpha\|_{\mathcal B^{\alpha+1}}
\geq\int_0^1\frac{(1-t)^b}{(1-t^2)^\alpha}\,dt
+(b+1)\sup_{0\leq r<1}(1-r^2)^{\alpha+1}
\int_0^1\frac{t(1-t)^b}{(1-rt)^{b+2}(1-t^2)^\alpha}\,dt.
$}
\end{equation}

Now we consider the upper bound. Suppose that $f\in H_\alpha^\infty$.
From the definition of the Korenblum norm,
\[
|f(t)|\leq\frac{\|f\|_{H_\alpha^\infty}}{(1-t^2)^\alpha},
\qquad 0\leq t<1.
\]
For $z\in\mathbb D$, the inequalities $|1-tz|\geq1-t|z|$
and \eqref{eq:Hb-Korenblum-Bloch-derivative} give
\[
\|\mathcal H_bf\|_{\mathcal B^{\alpha+1}}
=|\mathcal H_bf(0)|+\sup_{z\in\mathbb D}(1-|z|^2)^{\alpha+1}|(\mathcal H_bf)'(z)|,
\]
\[
\resizebox{0.94\linewidth}{!}{$\displaystyle
\|\mathcal H_bf\|_{\mathcal B^{\alpha+1}}
\leq\biggl[\int_0^1\frac{(1-t)^b}{(1-t^2)^\alpha}\,dt
+(b+1)\sup_{0\leq r<1}(1-r^2)^{\alpha+1}
\int_0^1\frac{t(1-t)^b}{(1-rt)^{b+2}(1-t^2)^\alpha}\,dt\biggr]
\|f\|_{H_\alpha^\infty}.
$}
\]
Together with \eqref{eq:Hb-Korenblum-Bloch-lower}, this proves
\eqref{eq:Hb-Korenblum-Bloch-norm}, provided that the supremum is
finite. The same inequalities show that $f_\alpha$ attains the
operator norm.

We next determine the boundary contribution of the second term.
Changing variables by $t=1-(1-r)x$, we obtain
\begin{equation}\label{eq:Hb-Korenblum-Bloch-scaled}
\begin{aligned}
&(b+1)(1-r^2)^{\alpha+1}
 \int_0^1\frac{t(1-t)^b}{(1-rt)^{b+2}(1-t^2)^\alpha}\,dt\\
&\quad=(b+1)(1+r)^{\alpha+1}
\int_0^{1/(1-r)}
\frac{[1-(1-r)x]x^{b-\alpha}}
{[2-(1-r)x]^\alpha(1+rx)^{b+2}}\,dx.
\end{aligned}
\end{equation}
For $\frac12\leq r<1$, the integrand on the right, extended by
zero for $x>1/(1-r)$, is dominated by
\[
\frac{x^{b-\alpha}}{(1+x/2)^{b+2}}.
\]
This function is integrable on $(0,\infty)$, since
$b-\alpha>-1$ and $\alpha>0$. Thus dominated convergence gives
\begin{equation}\label{eq:Hb-Korenblum-Bloch-boundary}
\begin{aligned}
&\lim_{r\to1^-}(b+1)(1-r^2)^{\alpha+1}
 \int_0^1\frac{t(1-t)^b}{(1-rt)^{b+2}(1-t^2)^\alpha}\,dt\\
&\quad=2(b+1)\int_0^\infty\frac{x^{b-\alpha}}{(1+x)^{b+2}}\,dx\\
&\quad=2(b+1)B(b+1-\alpha,1+\alpha).
\end{aligned}
\end{equation}
In particular, the right side of \eqref{eq:Hb-Korenblum-Bloch-norm}
is at least the first term in \eqref{eq:Hb-Korenblum-Bloch-bounds}.

We also verify its finiteness for every $0<\alpha<b+1$.
For $0\leq r\leq\frac12$, the second integral in
\eqref{eq:Hb-Korenblum-Bloch-norm}, including its prefactor, is
bounded by
\[
(b+1)2^{b+2}\int_0^1\frac{t(1-t)^b}{(1-t^2)^\alpha}\,dt<\infty.
\]
For $\frac12\leq r<1$, the preceding integrable majorant and
\eqref{eq:Hb-Korenblum-Bloch-scaled} give a bound independent of $r$.
This proves the boundedness assertion for $0<\alpha<b+1$.

To establish the exact norm for $0<\alpha\leq\frac23$, we now make
the change of variables
\[
x=\frac{1-t}{1-rt},\qquad t=\frac{1-x}{1-rx}
\]
in the second integral of \eqref{eq:Hb-Korenblum-Bloch-norm}. After
simplification, we obtain
\begin{equation}\label{eq:Hb-Korenblum-Bloch-substitution}
\begin{aligned}
&(b+1)(1-r^2)^{\alpha+1}
\int_0^1\frac{t(1-t)^b}{(1-rt)^{b+2}(1-t^2)^\alpha}\,dt\\
&\quad=(b+1)(1+r)^{\alpha+1}
\int_0^1\frac{(1-t)t^{b-\alpha}(1-rt)^{2\alpha-1}}
{[2-(1+r)t]^\alpha}\,dt.
\end{aligned}
\end{equation}
For $0<\alpha\leq\frac23$, let
\[
g(r,t)=\frac{(1-rt)^{2\alpha-1}}{[2-(1+r)t]^\alpha}.
\]
A direct differentiation yields
\[
\frac{\partial g(r,t)}{\partial r}
=\frac{t(1-rt)^{2\alpha-2}}{[2-(1+r)t]^{\alpha+1}}
\bigl[(1-\alpha)(1-rt)+(1-2\alpha)(1-t)\bigr].
\]
When $0<\alpha\leq\frac12$, both terms in the brackets are
nonnegative. When $\frac12<\alpha\leq\frac23$, we use
$1-rt\geq1-t$ to obtain
\[
(1-\alpha)(1-rt)+(1-2\alpha)(1-t)
\geq(2-3\alpha)(1-t)\geq0.
\]
Thus $g(r,t)$ is nondecreasing in $r$ for every $0<t<1$.
Since $(1+r)^{\alpha+1}$ is strictly increasing, the right side
of \eqref{eq:Hb-Korenblum-Bloch-substitution} is strictly increasing
in $r\in[0,1)$. Hence its supremum is the limit as $r$ tends to $1$
from below. By \eqref{eq:Hb-Korenblum-Bloch-boundary}, this supremum
equals $2(b+1)B(b+1-\alpha,1+\alpha)$. Formula
\eqref{eq:Hb-Korenblum-Bloch-sharp} now follows.

Finally, suppose $\frac23<\alpha<\min\{2,b+1\}$.
Since $2\alpha-1>0$, we have
\[
(1-rt)^{2\alpha-1}\leq1,\qquad
2-(1+r)t\geq2(1-t),\qquad 1+r\leq2.
\]
Using these inequalities in
\eqref{eq:Hb-Korenblum-Bloch-substitution}, we obtain
\[
\begin{aligned}
&(b+1)(1-r^2)^{\alpha+1}
\int_0^1\frac{t(1-t)^b}{(1-rt)^{b+2}(1-t^2)^\alpha}\,dt\\
&\quad\leq 2(b+1)\int_0^1t^{b-\alpha}(1-t)^{1-\alpha}\,dt\\
&\quad=2(b+1)B(b+1-\alpha,2-\alpha).
\end{aligned}
\]
Combining this estimate with the boundary lower bound
\eqref{eq:Hb-Korenblum-Bloch-boundary} proves
\eqref{eq:Hb-Korenblum-Bloch-bounds}.

If $\alpha\geq b+1$, the function
$f_\alpha(z)=(1-z^2)^{-\alpha}$ still has unit norm in
$H_\alpha^\infty$, but the integral at the origin diverges:
\[
\mathcal H_bf_\alpha(0)
=\int_0^1(1-t)^{b-\alpha}(1+t)^{-\alpha}\,dt
=+\infty.
\]
Indeed, the integrand is comparable to $(1-t)^{b-\alpha}$
near $t=1$, which is not integrable when $b-\alpha\leq-1$.
Thus $\mathcal H_b$ cannot be a bounded operator from
$H_\alpha^\infty$ to $\mathcal B^{\alpha+1}$ in this range.
\end{proof}

\begin{remark}
Setting \(b=0\) in Theorems 6.1, we obtain Theorems 7.1  in \cite{HuYe2025}.
\end{remark}

\section[Norm on the logarithmically weighted Korenblum space]{Norm of $\mathcal H_b$ on $H_{\alpha,\log}^{\infty}$}

In this section, we determine  the norm of the generalized Hilbert matrix
operator $\mathcal H_b$ acting on the logarithmically weighted Korenblum
space $H_{\alpha,\log}^{\infty}$. Recall that the Beta function is defined by
\[
B(s,t)=\int_0^1 x^{s-1}(1-x)^{t-1}\,dx,
\qquad \operatorname{Re}s>0,\quad\operatorname{Re}t>0,
\]
and satisfies
\[
B(s,t)=\frac{\Gamma(s)\Gamma(t)}{\Gamma(s+t)},
\]
where $\Gamma$ denotes the Gamma function. We will also use Euler's
reflection formula
\[
\Gamma(z)\Gamma(1-z)=\frac{\pi}{\sin(\pi z)},
\qquad z\in\mathbb C\setminus\mathbb Z.
\]
In particular, for $b\geq0$ and $0<\alpha<1$,
\[
\int_0^1 t^{b-\alpha}(1-t)^{\alpha-1}\,dt
=B(1+b-\alpha,\alpha)
=\frac{\Gamma(1+b-\alpha)\Gamma(\alpha)}{\Gamma(1+b)}.
\]

\begin{theorem}\label{thm:Hb-log-self}
For $b>0$ and $0<\alpha<1$, the generalized Hilbert matrix operator
$\mathcal H_b$ is bounded on $H_{\alpha,\log}^{\infty}$.
Moreover,
\begin{equation}\label{eq:Hb-log-self-norm}
\|\mathcal H_b\|_{H_{\alpha,\log}^{\infty}\to H_{\alpha,\log}^{\infty}}
=\sup_{0\leq r<1}\int_0^1
\frac{(1-t)^b(1+r)^\alpha[1-(1-t)r]^{2\alpha-1}
\log\!\frac{2e^{1/\alpha}}{1-r^2}}
{(1-t)^\alpha[1+t-(1-t)r]^\alpha
\log\!\frac{2e^{1/\alpha}}{1-\phi_t(r)^2}}\,dt.
\end{equation}
Furthermore,
\begin{equation}\label{eq:Hb-log-self-beta-lower}
\|\mathcal H_b\|_{H_{\alpha,\log}^{\infty}\to H_{\alpha,\log}^{\infty}}
\geq B(1+b-\alpha,\alpha)
=\frac{\Gamma(1+b-\alpha)\Gamma(\alpha)}{\Gamma(1+b)}.
\end{equation}
\end{theorem}

\begin{proof}
Let
\[
f_\alpha(z)=\frac{1}{(1-z^2)^\alpha
\log\!\frac{2e^{1/\alpha}}{1-z^2}}.
\]
From the proof of Theorem~\ref{thm:Hb-radial-norm}, we know that
$\|f_\alpha\|_{H_{\alpha,\log}^{\infty}}=1$, and
\[
\begin{aligned}
T_tf_\alpha(z)
&=\frac{(1-t)^b}{1-(1-t)z}
 f_\alpha\!\left(\frac{t}{1-(1-t)z}\right)\\
&=\frac{(1-t)^b}{1-(1-t)z}
\frac{1}{\left(1-\left(\frac{t}{1-(1-t)z}\right)^2\right)^\alpha
\log\!\frac{2e^{1/\alpha}}
{1-\left(\frac{t}{1-(1-t)z}\right)^2}}.
\end{aligned}
\]
Since $\mathcal H_bf_\alpha(z)=\int_0^1T_tf_\alpha(z)\,dt$,
we obtain
\begin{align}
&\|\mathcal H_b\|_{H_{\alpha,\log}^{\infty}\to H_{\alpha,\log}^{\infty}}
\geq\|\mathcal H_bf_\alpha\|_{H_{\alpha,\log}^{\infty}}
\notag\\
&=\sup_{z\in\mathbb D}(1-|z|^2)^\alpha
\log\frac{2e^{1/\alpha}}{1-|z|^2}
\left|\int_0^1T_tf_\alpha(z)\,dt\right|\notag\\
&\geq\sup_{0\leq r<1}(1-r^2)^\alpha
\log\frac{2e^{1/\alpha}}{1-r^2}
\int_0^1\frac{(1-t)^b}{1-(1-t)r}
\frac{dt}{\left(1-\left(\frac{t}{1-(1-t)r}\right)^2\right)^\alpha
\log\!\frac{2e^{1/\alpha}}
{1-\left(\frac{t}{1-(1-t)r}\right)^2}}\notag\\
&=\sup_{0\leq r<1}\int_0^1
\frac{(1-t)^b(1+r)^\alpha
\bigl((t-1)r+1\bigr)^{2\alpha-1}
\log\!\frac{2e^{1/\alpha}}{1-r^2}}
{(1-t)^\alpha\bigl((t-1)r+1+t\bigr)^\alpha
\log\!\frac{2e^{1/\alpha}}
{1-\left(\frac{t}{1-(1-t)r}\right)^2}}\,dt.
\label{eq:Hb-log-self-lower}
\end{align}
Here, on the positive radius all factors in the integral are positive,
and we used
\[
1-\left(\frac{t}{1-(1-t)r}\right)^2
=\frac{(1-r)(1-t)\bigl((t-1)r+1+t\bigr)}
{\bigl((t-1)r+1\bigr)^2}.
\]

Next, we estimate the upper bound. By \eqref{eq:Hb-T-upper},
\[
|T_tf(z)|\leq\frac{|w_t(z)|}{
(1-|\phi_t(z)|^2)^\alpha
\log\!\frac{2e^{1/\alpha}}{1-|\phi_t(z)|^2}}
\|f\|_{H_{\alpha,\log}^{\infty}}.
\]
Since $|w_t(z)|\leq w_t(|z|)$,
$|\phi_t(z)|\leq\phi_t(|z|)$, and
$x^\alpha\log(2e^{1/\alpha}/x)$ is increasing for $0<x<2$,
the preceding estimate gives
\[
|T_tf(z)|\leq\frac{w_t(r)}{
(1-\phi_t(r)^2)^\alpha
\log\!\frac{2e^{1/\alpha}}{1-\phi_t(r)^2}}
\|f\|_{H_{\alpha,\log}^{\infty}},
\qquad r=|z|.
\]
Thus
\begin{align}
&\|\mathcal H_bf\|_{H_{\alpha,\log}^{\infty}}
=\sup_{z\in\mathbb D}(1-|z|^2)^\alpha
\log\frac{2e^{1/\alpha}}{1-|z|^2}
\left|\int_0^1T_tf(z)\,dt\right|\notag\\
&\leq\sup_{z\in\mathbb D}(1-|z|^2)^\alpha
\log\frac{2e^{1/\alpha}}{1-|z|^2}
\int_0^1|T_tf(z)|\,dt\notag\\
&\leq\sup_{0\leq r<1}(1-r^2)^\alpha
\log\frac{2e^{1/\alpha}}{1-r^2}
\int_0^1\frac{(1-t)^b}{1-(1-t)r}
\frac{dt}{(1-\phi_t(r)^2)^\alpha
\log\!\frac{2e^{1/\alpha}}{1-\phi_t(r)^2}}
\|f\|_{H_{\alpha,\log}^{\infty}}\notag\\
&=\sup_{0\leq r<1}\int_0^1
\frac{(1-t)^b(1+r)^\alpha
\bigl((t-1)r+1\bigr)^{2\alpha-1}
\log\!\frac{2e^{1/\alpha}}{1-r^2}}
{(1-t)^\alpha\bigl((t-1)r+1+t\bigr)^\alpha
\log\!\frac{2e^{1/\alpha}}
{1-\left(\frac{t}{1-(1-t)r}\right)^2}}\,dt\;\|f\|_{H_{\alpha,\log}^{\infty}}.
\label{eq:Hb-log-self-upper}
\end{align}
Combining \eqref{eq:Hb-log-self-lower} and
\eqref{eq:Hb-log-self-upper} gives \eqref{eq:Hb-log-self-norm},
provided that the displayed supremum is finite.

We now verify finiteness and obtain the lower bound. Using the
integral representation \eqref{eq:Hb-integral}, we have
\[
\begin{aligned}
&(1-r^2)^\alpha\log\frac{2e^{1/\alpha}}{1-r^2}
\mathcal H_bf_\alpha(r)\\
&\quad=(1-r^2)^\alpha\log\frac{2e^{1/\alpha}}{1-r^2}
\int_0^1\frac{(1-t)^b}{(1-rt)^{b+1}}
\frac{dt}{(1-t^2)^\alpha
\log\!\frac{2e^{1/\alpha}}{1-t^2}}.
\end{aligned}
\]
For $r\geq\frac12$, substitute
\[
u=\frac{1-t}{1-r},\qquad t=1-(1-r)u.
\]
Then
\begin{equation}\label{eq:Hb-boundary-integral}
\begin{aligned}
&(1-r^2)^\alpha\log\frac{2e^{1/\alpha}}{1-r^2}
\mathcal H_bf_\alpha(r)\\
&=\int_0^{1/(1-r)}
\frac{(1+r)^\alpha u^{b-\alpha}}
{(1+ru)^{b+1}[2-(1-r)u]^\alpha}
\frac{\log\!\frac{2e^{1/\alpha}}{(1-r)(1+r)}}
{\log\!\frac{2e^{1/\alpha}}
{(1-r)u[2-(1-r)u]}}\,du.
\end{aligned}
\end{equation}
For each fixed $u>0$, the integrand in
\eqref{eq:Hb-boundary-integral} tends, as $r$ tends to $1$
from below, to
\[
\frac{u^{b-\alpha}}{(1+u)^{b+1}}.
\]
To justify passage to the limit, note that if $0<u\leq1$, then
$t\geq r$, so
\[
0<\frac{\log\!\frac{2e^{1/\alpha}}{1-r^2}}
{\log\!\frac{2e^{1/\alpha}}{1-t^2}}\leq1.
\]
If $u>1$, then $t<r$ and
\[
\begin{aligned}
\frac{\log\!\frac{2e^{1/\alpha}}{1-r^2}}
{\log\!\frac{2e^{1/\alpha}}{1-t^2}}
&=1+\frac{\log\!\frac{1-t^2}{1-r^2}}
{\log\!\frac{2e^{1/\alpha}}{1-t^2}}\\
&\leq1+\frac{\log(2u)}{\log(2e^{1/\alpha})}.
\end{aligned}
\]
For $r\geq\frac12$ we also have
$1+ru\geq1+u/2$, $1+r\leq2$, and
$2-(1-r)u=1+t\geq1$.
Consequently, the integrand in
\eqref{eq:Hb-boundary-integral}, extended by zero outside its
integration interval, is bounded by
\[
\frac{2^\alpha u^{b-\alpha}}{(1+u/2)^{b+1}}
\left(1+\frac{\log(2(1+u))}{\log(2e^{1/\alpha})}\right).
\]
This bound is integrable on $(0,\infty)$: near $0$ it behaves
like $u^{b-\alpha}$ with $b-\alpha>-1$, while at infinity it is
$O(u^{-1-\alpha}\log(2+u))$ with $\alpha>0$.
The dominated convergence theorem therefore gives
\begin{align}
\lim_{r\to1^-}(1-r^2)^\alpha
\log\frac{2e^{1/\alpha}}{1-r^2}\mathcal H_bf_\alpha(r)
&=\int_0^\infty\frac{u^{b-\alpha}}{(1+u)^{b+1}}\,du\notag\\
&=B(1+b-\alpha,\alpha).
\label{eq:Hb-boundary-beta}
\end{align}
The same domination shows that the expression on the left of
\eqref{eq:Hb-boundary-integral} is uniformly bounded for
$\frac12\leq r<1$.
For $0\leq r\leq\frac12$, the factor
$1/(1-rt)^{b+1}$ is bounded by $2^{b+1}$, and the remaining
integral is finite because $b-\alpha>-1$.
Thus the supremum in \eqref{eq:Hb-log-self-norm} is finite.
Finally, \eqref{eq:Hb-log-self-lower} and
\eqref{eq:Hb-boundary-beta} imply
\[
\|\mathcal H_b\|_{H_{\alpha,\log}^{\infty}\to H_{\alpha,\log}^{\infty}}
\geq B(1+b-\alpha,\alpha),
\]
which proves \eqref{eq:Hb-log-self-beta-lower} and completes the
proof.
\end{proof}
\begin{remark}
Setting \(b=0\) in Theorems 7.1, we obtain Theorems 4.1  in \cite{HuYe2025}.
\end{remark}

\subsection*{Conflicts of Interest}
The authors declare that there are no conflicts of interest regarding the publication of this paper.
\subsection*{Funding Information}
The second author was supported by    National Natural Science Foundation of China (Grant No. 11671357) and Zhejiang Province Natural Science Foundation of China(Grant No. LY23A010003).
\subsection*{Author contribution}
All authors have contributed equally to the writing of this paper.
\subsection*{Data Availability Statement}

The authors declare that all data and material in this paper are available.

\end{document}